\documentclass[12pt]{article}
\usepackage{amssymb}
\usepackage{amsmath}
\usepackage{theorem}
\input xy
\xyoption{all}

\newtheorem{thm}{Theorem}[section]
\newtheorem{prop}[thm]{Proposition}
\newtheorem{prob}[thm]{Problem}
\newtheorem{coro}[thm]{Corollary}
\newtheorem{lemm}[thm]{Lemma}

{\theorembodyfont{\rm}

\newtheorem{exam}[thm]{Example}
\newtheorem{rem}[thm]{Remark}

}

\def\ra{\rightarrow}

\def\C{{\mathbb C}}

\def\P{{\mathbb P}}
\def\Q{{\mathbb Q}}

\def\Z{{\mathbb Z}}
\def\C{{\mathbb C}}

\def\Alb{{\rm Alb}}
\def\alb{{\rm alb}}
\def\Pic{{\rm Pic}}
\def\E{{\cal E}}
\def\SL{{\rm SL}}

\def\Sym{{\rm Sym}}
\def\Ex{{\Delta}}
\def\Jac{{\mathcal J}}
\def\rk{{\rm rk}}
\def\W{{\mathcal W}}
\def\Proj{{\rm Proj}}
\def\NS{{\rm NS}}

\author{Fedor Bogomolov\\
\small  Courant Institute of Mathematical Sciences, N.Y.U. \\
\small 251 Mercer str. \\
\small New York, NY 10012, U.S.A.\\
\small e-mail: bogomolo@cims.nyu.edu\\
\small  and\\
Yuri Tschinkel\\
\small Department of Mathematics, \\
\small Princeton University\\
\small Fine Hall, Washington Road\\
\small Princeton, NJ 08544-1000,  U.S.A.  \\
\small e-mail: ytschink@math.princeton.edu
}

\title{Lagrangian subvarieties of abelian fourfolds}

\begin{document} 
 
\date{\today}

\maketitle

\pagebreak

\centerline{\it Dedicated to Professor K. Kodaira}

\

\section{Introduction}

\

\

Let $({\mathcal W}, \omega)$ be a smooth projective algebraic
variety of dimension $2n$ over $\C$ together with a holomorphic
$(2,0)$-form of maximal rank $2n$. A subvariety $X\subset {\mathcal W}$
is called {\em weakly lagrangian} if $\dim X \le n$ and if 
the restriction of $\omega$ to $X$ is trivial (notice that
$X$ can be singular). 
An $n$-dimensional subvariety 
$X\subset \W$ with this property is called {\em lagrangian}. 
For example, any curve $C$ contained in 
a K3 or abelian surface $S$ is lagrangian.  
Further examples of lagrangian subvarieties are
obtained by taking a curve $C\subset S$ and by
considering the corresponding symmetric products. 
Alternatively, one could look at a product of 
different curves (of genus $> 1$) inside a product
of abelian varieties. 
We will say that a variety
$X\subset \W$ is {\em fibered} if it admits a dominant map
onto a curve of genus $ > 1$. 
In this note we construct examples of nonfibered 
lagrangian surfaces in abelian varieties. 
Our motivation comes from the following

\begin{prob}
Find examples of projective surfaces
with a nontrivial fundamental group. In particular, 
find examples where the fundamental group
has a nontrivial nilpotent tower.  
\end{prob}

If $X$ is fibered over a curve 
$C$ of genus $>1$ then the fundamental group 
$\pi_1(X)$ surjects onto a subgroup of
finite index in $\pi_1(C)$ and consequently
both $\pi_1(X)$ and its nilpotent tower are big.
Therefore, we are interested in examples of surfaces where the 
nontriviality of $\pi_1(X)$ is not induced from curves. 
Consider the map to the Albanese variety $\alb\,:\, X\ra \Alb(X)$.  
One is interested in situations where the natural map
$$
H^2(\Alb(X),\C)\ra H^2(X,\C)
$$
has a nontrivial kernel -- the triviality of the  kernel implies  
the triviality of the nilpotent tower (tensor $\Q$).  
Such examples were given by  
Campana (\cite{campana}, Cor. 1.2) and 
Sommese-Van de Ven (\cite{sommese-vandeven}). 
However, there the kernel was found in the map
$$
\Pic(\Alb(X))\ra \Pic(X).
$$
In their construction the fundamental group of the 
variety is a central extension of an abelian group 
(and the lower central series has only two steps).
The lagrangian property of $\alb(X)\subset \Alb(X)$ implies
that there is a nontrivial kernel in $H^{2,0}$ (rather than 
on the level of the Picard groups).  
We produce an infinite series of surfaces $X$ 
of different topological types which 
are contained in abelian varieties and 
are lagrangian (or weakly lagrangian)
with respect to a nondegenerate $(2,0)$-form. 
The weakly lagrangian property for a subvariety $X$
of an abelian variety $A$ is related to the nontriviality 
of the fundamental group $\pi_1(X)$ as follows:
the number of linearly independent generators 
of the second quotient of the central series 
of $\pi_1(X)$ is bounded from below
by the dimension of the rational envelope of
the space of those $(2,0)$-forms on $A$ which restrict
trivially to $X$.

Our construction uses dominant maps 
between K3 surfaces. Let $S$ be a K3 surface 
and $g_1,g_2$ dominant rational maps of $S$ to 
Kummer K3 surfaces $S_1,S_2$ (blowups of quotients of
abelian surfaces $A_1,A_2$ by standard involutions). 
Then the (birational) preimage $X$ of 
$(g_1,g_2)(S)\subset S_1\times S_2$
in $A_1\times A_2$ is a lagrangian 
surface (in general, singular). For special choices of 
$g_1,g_2$ we can compute some basic invariants of $X$
and, in particular, show that $X$ is nonfibered.
For example, let $A$ be an abelian surface which is not 
isogenous to a product of elliptic curves and $S$ 
the associated Kummer surface.  Assume that
$g_1\,:\,S\ra S$ is an isomorphism 
and $g_2\,:\,S\ra S$ is not an isomorphism. 
Then $X\subset A\times A$ is not fibered, $\Alb(X)$ is isogenous
to $A\times A$ and $X$ is lagrangian with 
respect to exactly one 2-form on $A\times A$.  
We analyze other constructions, with $S$ an elliptic 
Kummer surface and $g_j$ induced from the group law.
We don't determine the actual structure of $\pi_1(X)$ 
(it presumably depends on $X$), but it seems quite 
plausible that for some $X$ from our list 
$\pi_1(X)$ has a rather nontrivial nilpotent tower.

\

{\bf Acknowledgments.} The first author was partially
supported by the NSF. The second author was partially supported
by the NSA.  The paper was motivated by a
question raised by F.~Catanese. We would like
to thank him for useful discussions in the 
early stages of this work. We thank the referee for comments
which helped to improve the exposition.

\

\section{Preliminaries}
\label{sect:2}

Let $V_{\Q}$ be a finite dimensional $\Q$-vector space and 
$V_{\C}$ its complexification. Let $w\in V_{\C}$ be 
a vector. We denote by 
$L(w)\subset V_{\Q}$ its rational envelope, i.e.,
the smallest linear subspace such that $w\in L(w)_{\C}$.  
More generally, if $W\subset V_{\C}$ is any
set, we will denote by $L(W)\subset V_{\Q}$ the smallest
linear subspace such that $L(w)\subset L(W)$ for all $w\in W$. 
An element $w\in V_{\C}$ will be
called {\em k-generic} if the dimension of $L(w)$
is $\dim V_{\Q} - k$.

\begin{rem}
The set of rationally defined subspaces in $V_{\C}$ is 
countable. Therefore, for any linear subspace $W\subset V_{\C}$
we have $L(W)=L(w)$ for all  $w\in W$ which 
are not contained in a countable number of linear subspaces.
\end{rem}

\begin{prop}
Let $(A,\omega)$ be an abelian variety of dimension $2n$
together with a nondegenerate holomorphic $(2,0)$-form $\omega$. 
Assume that $\omega$ defines a $k$-generic 
element in $H^{2,0}(A,\C)$ with $k  < 3 $. Then there are no 
weakly lagrangian
surfaces in $A$. 
\end{prop}

{\em Proof.}
Assume that we have a weakly lagrangian surface 
$i\,:\, X\hookrightarrow A$. Consider the induced 
$\Q$-rational homomorphism
$$
i^* \,:\, H^2(A,\Q)\ra H^2(X,\Q).
$$
The class of $\omega$ is contained in 
the ($\Q$-rationally defined)
kernel of $i^*$. Since $\omega$ is $k$-generic in $H^{2,0}(A,\C)$
for $k < 3$ the image of $H^2(A,\C) $ 
in $H^2(X,\C)$ has rank at most
$2$. If we had  a  
holomorphic form $w'\in H^{2, 0}(A,\C) $ 
restricting nontrivially to $X$ then the image of $H^2(A,\C) $ in
$H^2(X,\C) $ would contain at least three linearly independent
forms $ w', \bar w'$ and the image of 
a polarization from $\Pic(A)$.
Thus the triviality of $\omega$ on $X$ implies that
{\em all} $(2,0)$-forms restrict trivially to $X$. 
Consequently, the dimension of $X$ is $\le 1$, contradiction.   
{}\hfill $\square$

\

We end this section with a simple description
of maps of fibered surfaces onto curves (see also \cite{catanese}). 

\begin{prop}\label{lemm:factoring} 
For every smooth projective 
algebraic surface $V$ there exists
a universal algebraic variety ${\cal U}(V)$ such that
every dominant rational map $V\ra C$ onto a smooth
curve of genus $\ge 2$ factors through ${\cal U}(V)$. 
\end{prop}

{\em Proof.} 
First observe that the set of dominant rational maps
of $V$ onto 
smooth curves with generically irreducible fibers is
countable. Indeed, the class $f^*(L)$ (where $L$ is 
some polarization on the image-curve) defines $f$.    
If two such classes $[f]$
and $[f']$ differ by an element in $\Pic^0(V)$, 
then the maps must be the same. Otherwise, 
some fiber of $f$ would surject onto the image of $f'$
and therefore, the degrees of the two classes
$[f]$  and $[f']$ on this
fiber would be  different.

Further, dominant maps onto smooth curves of genus $\ge 2$
define a linear subspace $W\subset H^0(V,\Omega^1)$ by
the property that wedge-products of linearly independent
forms $w,w'\in W$ are trivial. 
The set of maximal subspaces $W$ in $H^0(V,\Omega^1)$
with this property is an algebraic variety. 
Remark that for two forms $\omega, \omega'$
with $\omega\wedge \omega'=0$  their ratio $\omega/\omega'$ is 
a nonconstant rational function, which is constant on the fibers
of the foliations on $V$ defined by $\omega$, resp. $\omega'$. 
Thus, any such subspace $W$ defines a foliation ${\cal F}_W$ 
on $V$ with compact nonintersecting fibers (locally the 
form defines
a holomorphic map, its fibers are the fibers of ${\cal F}_W$).
Therefore, ${\cal F}_W$ defines a dominant rational map 
onto a curve $C_W$ (which lifts to a map of $V$ onto the
normalization of $C_W$).
Different $W$ define different foliations
and different morphisms
(the sum of two spaces $W$ and $W'$ with
the same foliations ${\cal F}_W={\cal F}_{W'}$ 
has the same property, contradicting the maximality of
$W$ and $W'$). Since the 
set of such spaces $W$ is on the one hand 
algebraic (a finite union of subvarieties of a Grassmannian) 
and on the other hand countable, it must be finite. \hfill $\square$ 
 
\begin{coro}
\label{coro:factoring} 
For every smooth projective 
algebraic surface $V$ there exists
a universal algebraic variety ${\cal U}(V)$ such that
every dominant rational map $V\ra V'$, where $V'$ is a 
product of smooth curves of genus $\ge 2$, 
factors through ${\cal U}(V)$. 
\end{coro}

{\em Proof.} We have shown in Lemma~\ref{lemm:factoring}
that there is a finite
number of dominant rational maps onto curves of genus $\ge 2$. 
The product of these curves and maps is ${\cal U}(V)$. 
Universality follows. \hfill $\square$

\section{Construction}
\label{sect:construction}

Let $A$ be an abelian surface and $\tau$ the involution 
$\tau(a) = -a$ (for $a\in A$). This involution 
acts on $A$ with 16 fixed points. 
We denote by $\tilde{S}$ the (singular) quotient $A/\tau$.
Let $A^*$ be the blowup of $A$ in the 16 points. 
The involution $\tau$ extends to a fixed point free 
action on $A^*$ and the quotient $A^*/\tau$ is a K3 (Kummer) surface $S$, 
a blowup of $\tilde{S}$. 
It contains 16 exceptional curves and we
will denote by $\Ex=\cup_{k=1}^{16} \Ex^k$ their union.
We shall denote by  $\delta\,:\, A\ra \tilde{S}$  the 
double cover. 
Every K3 surface has a unique (up to constants)
nondegenerate holomorphic $(2,0)$-form. 

\

We start with two {\em simple} abelian surfaces $A_1$ and 
$A_2$ (in particular, 
$A_j$ don't contain elliptic curves - this assumption simplifies
the discussion in Section~\ref{sect:details}). 
We consider the corresponding Kummer surfaces $S_j$ and 
we assume that there exists another K3 surface
$S$ together with two dominant rational maps $g_j\,:\, S\ra S_j$.
We denote by $\tilde{g}_j$ the induced maps $S\ra \tilde{S}_j$. 
Now consider the map 
$$
\tilde{\varphi}
=(\tilde{g}_1,\tilde{g}_2)\,:\, S\ra \tilde{S}_1\times \tilde{S}_2.
$$
We have two double covers 
$$
(\delta_1,1)\,:\,
A_1\times \tilde{S}_2 \ra \tilde{S}_1\times \tilde{S}_2
$$
and 
$$
(1,\delta_2)\,:\, \tilde{S}_1\times A_2\ra \tilde{S}_1\times \tilde{S}_2.
$$
Denote by $\tilde{Y}_1=(\delta_1,1)^{-1}(\tilde{\varphi}(S))$
and by $\tilde{Y}_2=(1,\delta_2)^{-1}(\tilde{\varphi}(S))$.
Let 
$$
\tilde{X}=\tilde{X}_{\tilde{\varphi}}:=(\delta_1,\delta_2)^{-1}(
\tilde{\varphi}(S))
$$ 
be the  preimage of 
$S$ in $A_1\times A_2$.
We see that the 
$(2,2)$-covering $\tilde{X}\ra \tilde{\varphi}(S)$
is 
$$
\tilde{X}=\tilde{Y}_1\times_{\tilde{\varphi}(S)}
\tilde{Y}_2\subset A_1\times A_2.
$$
The surface $\tilde{X}$ is, in general, singular. 

\begin{lemm}\label{lemm:largange}
There exists a nondegenerate holomorphic 2-form  $\omega'$ 
on $A_1\times A_2$ such that 
$\tilde{X}$ is lagrangian with respect to $\omega'$. 
\end{lemm}

{\em Proof.}
The singular surface $\tilde{S}_j$ carries a
nondegenerate holomorphic 2-form, which we again denote by 
$\omega_j$.
Evidently, $\tilde{g}_j^*\omega_j=\lambda_j\omega$ for 
some nonzero numbers $\lambda_j$. 
Consider the form 
$$
\omega'= \lambda_2 \omega_1 - 
\lambda_1 \omega_2
$$
on $\tilde{S}_1\times \tilde{S}_2$. The form $\omega'$ is
identically zero on  the image $\tilde{\varphi}(S)$
(since it is trivial on the open part of $S$ where $\tilde{\varphi}$
is smooth). 
Since both forms $\omega_j$ lift to nondegenerate forms on the 
abelian surfaces $A_j$ the form $\omega'$ lifts 
to a nondegenerate form on $A_1\times A_2$.
The restriction of the lift of $\omega'$ to $\tilde{X}$ is
identically zero on the smooth points of $\tilde{X}$. 
\hfill $\square$

\begin{exam} If $A_1=A_2=A$ and the maps 
$\tilde{g}_j={\rm id}$ then the resulting surface $\tilde{X}$ is
the abelian surface $A$, embedded diagonally into $A\times A$. 
\end{exam}

\begin{exam}\label{exam:unique}
Let $S$ be a K3 surface, which is simultaneously Kummer and 
a double cover of $\P^2$. Denote 
by $\theta$ the covering involution on $S$. 
Put $S_1=S_2=S$, $\tilde{g}_1={\rm id}$ and 
$\tilde{g}_2=\theta$. 
Then the corresponding surface $\tilde{X}$ is lagrangian. A special case
of this construction is obtained as follows:
Let $C$ be a curve of genus $2$ 
and $\sigma$ a hyperelliptic involution on $C$.
Consider $C\times C$, together with the involutions:
$\sigma_{12}$ interchanging the factors and 
$$
\sigma_a\,:\,
(c_1,c_2)\ra (\sigma(c_2),\sigma(c_1)).
$$   
Denote by  $Y_1=C\times C/\sigma_{12}$ and $Y_2=C\times C/\sigma_a$. 
We see that both $Y_1$ and $Y_2$ are isomorphic to a symmetric square
of $C$, which is birational to the same abelian surface $A$. 
Denote by $\tilde{S}$ the quotient of $A$ by the standard involution $\tau$. 
Observe that $S$ is realized as a double cover of 
$\P^2=\P^1\times \P^1/\sigma_{12}$. Then 
$$
(\tau,\tau)^{-1}((1,\theta)(S))=C\times C\subset A\times A.
$$
\end{exam}

In the following sections we will show that 
we can arrange a situation where $\tilde{X}$ is not contained in 
any abelian subvariety of $A_1\times A_2$ and 
where it does not admit any dominant
morphisms onto a curve of genus $\ge 2$.

\section{Details and proofs}
\label{sect:details}

\subsection{Desingularization}
\label{sect:desingularization}

To analyze the surface $\tilde{X}$ 
constructed in Section~\ref{sect:construction}
we will need its sufficiently explicit partial 
desingularization.  It will be constructed in two steps. 

\

Let $\Sigma $ be the minimal finite set of points in $S$ such
that both maps $g_1,g_2$ are well defined on
the complement $S^0=S\setminus \Sigma$. 
On $S^0$ the maps $g_1$ and $g_2$ are local isomorphisms. 
Consider the divisors 
$\Ex_j\subset S_j$. 
The preimages $D_j^0$ of these divisors 
(under the maps $g_j$) in 
$S^0$ are smooth and the components of $D^0_j$ 
don't intersect in $S^0$. However, the divisors $D_1^0$ and $D^0_2$ do
intersect in $S^0$, and we denote by $p_1,...,p_l\in S^0$
their intersection points. Let $D_j$ be the closure
of $D^0_j$ in $S$. 
We fix a blowup  $\bar{S}$ of $S$  
with  centers supported in $\Sigma$ 
such that the preimage of the intersection 
$D_1\cap D_2$ in the neighborhood
of every point in $\Sigma$ is a normal crossing divisor.
Of course, $\bar{S}$ has a map to $\tilde{S}_j$ and we denote by
$\bar{Y}_j$ the fibered product $\bar{S}\times_{\tilde{S}_j}A_j$. 
The surface $\bar{Y}_j$ is a double cover
of $\bar{S}$ and it has at most $A1$-singularities (double points), since
$D_j^0$ is smooth in $S^0$. Consider the fibered product 
$\bar{X}= \bar{Y}_1\times_{\bar{S}}\bar{Y}_2$. All singular
points of $\bar{X}$ which are not of type $A_1$ lie
over the intersection points of $D_1^0$ and $D_2^0$ in $S^0$.
The surface $\bar{X}$ has a natural action of $\Z/2+\Z/2$. It admits 
equivariant surjective maps $\bar{\delta}_j\,:\, \bar{X}\ra\bar{Y}_j$, 
where $\bar{Y}_j$ is the quotient of $\bar{X}$ by 
the involution $\bar{\tau}_j$. We denote by $\bar{Y}_{12}$ the
quotient of $\bar{X}$ by 
$\bar{\tau}_{12}=\bar{{\tau}}_1\bar{{\tau}}_2$, 
it is still singular. 

\

Let $\hat{S}$ be the minimal 
blowup of $\bar{S}$ with support 
in the points $p_1,...,p_l\in S^0\subset \bar{S}$
such that proper transforms of the irreducible
components $D_1^0$ and $D_2^0$ are disjoint in the preimage of $S^0$ in 
$\hat{S}$. Now we define $Y_j$ as 
the induced (from the open part) double covers of $\hat{S}$. 
Their ramification is contained in 
the full transforms of the divisors
$D_j$. Define $X$ as the fibered product 
$Y_1\times_{\hat{S}}Y_2$. We have
the induced involutions 
(again denoted by  $\tau_j$, $\tau_{12}$) on $X$ and 
we define
$Y_{12} $ as the quotient of $X$ under $\tau_{12}$
(it admits a map onto $\hat{S}$). 
By construction, the 
surfaces $X$, $Y_j$ and $Y_{12}$ all have 
at most $A1$-singularities.  
One has surjective regular maps $X \ra \tilde{X},
Y_j\ra\tilde{Y}_j$ etc.

\

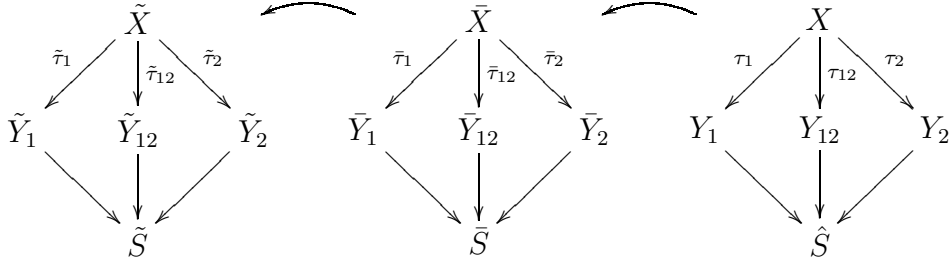
\centerline{
\xymatrix{
          &  \tilde{X}
\ar[dl]_{\tilde{\tau}_1}\ar[d]^{\tilde{\tau}_{12}}
\ar[dr]^{\tilde{\tau}_2} 
&    &  {}\ar@/_/[l]  &  \bar{X} 
\ar[dl]_{\bar{\tau}_1}\ar[d]^{\bar{\tau}_{12}}\ar[dr]^{\bar{\tau}_2} 
&    & {}\ar@/_/[l]   & X 
\ar[dl]_{{\tau}_1}\ar[d]^{{\tau_{12}}}\ar[dr]^{{\tau}_2}
&    
\\
  \tilde{Y}_1  \ar[dr]     &  \tilde{Y}_{12}  \ar[d] &   
\tilde{Y}_2\ar[dl] & 
   \bar{Y}_1 \ar[dr] &  \bar{Y}_{12}    \ar[d] &   \bar{Y}_2\ar[dl] &
       {Y}_1 \ar[dr] &  {Y}_{12}        \ar[d] &   {Y}_2   \ar[dl]    
\\
         &    \tilde{S}    &  &  & \bar{S} &  &  & \hat{S} 
}
}

\

\begin{lemm}\label{lemm:compon}
For all $k=1,...,l$ every irreducible component 
of the preimage of $p_k$ in $Y_{12}$
is a rational curve.
\end{lemm}

{\em Proof.}
Notice that $D_j^0$ are smooth in $S^0$ and that their
components don't intersect in $S^0$. Every point $p_k$
is a point of intersection of an irreducible component of $D_1^0$
with an irreducible component of $D_2^0$. 
In the neighborhood of $p_k$ these two divisors 
have a canonical form: $D_1$ is given by $x=0$ and $D_2$ is given 
by $y=x^n$, where $n$ is the order of tangency. 
There is a standard
chain of blowups separating the proper transforms of $D_1$ and $D_2$ 
over $p_k$ (which intersect one of the ends of the chain in two 
distinct points). The induced double cover on every component
of the preimage of $p_k$ is ramified in at most three points (hence
in fact, two) and is therefore a rational curve. 
\hfill $\square$

\begin{coro}\label{coro:factors}
The map $\alb\,:\, Y_{12}\ra \Alb(Y_{12})$ factors
through $\bar{Y}_{12}$.
\end{coro}

{\em Proof.} Indeed, the natural map 
$Y_{12}\ra \bar{Y}_{12}$ contracts  
connected graphs of rational curves to distinct points in 
$\bar{Y}_{12}$. 
Since these connected graphs of rational curves 
map into points in $\Alb(Y_{12})$  our claim follows. 
\hfill $\square$  

\

\subsection{Elliptic fibrations}
\label{sect:fibr}

Let $\E\ra \P^1$ be a Jacobian elliptic fibration and $M_1, M_2$
two irreducible 
horizontal divisors on $\E$. Let $M_{12}\subset \E$ 
be the divisor
of pairwise differences: $M_{12}\cap \E_b$ is the set of all points
of the form $p_1-p_2$, where $p_j\in  M_j\cap \E_b$.
The divisor $M_{12}$ may
have several irreducible components.
We shall say that the divisor $M_{12}$ is torsion
if every irreducible component of $M_{12}\subset \E$ 
consists of torsion points.   

\begin {lemm} 
\label{lemm:restriction}  
Consider the 
restriction map to the generic fiber 
$$
\eta^*\,:\, \Pic(\E)\ra \Pic(\E_{\eta}).
$$ 
If the divisor $M_{12}$ is torsion  
then there exists a positive integer 
$N$ such that the class  
$\eta^* ([M_{12}])\in \Pic^{(0)}(\E_{\eta})$ 
is annihilated by $N$.
\end{lemm} 

{\em Proof.}
For every irreducible component of $M_{12}$ there exists
a positive integer $N'$ such that all points $p$ in this component
are annihilated by $N'$. 
Since the divisor $M_{12}$ has only a finite 
number of irreducible components we can find an $N$
annihilating all points in $M_{12}$.
Consider the class $N[M_{12}]$. The
corresponding divisor is trivial upon restriction to 
the generic fiber $\E_{\eta}$. \hfill $\square$  

\begin{coro} \label{coro:eta}
The kernel of the map 
$\eta^*\,:\,\Pic(\E)\ra \Pic(\E_{\eta})$
is isomorphic to the subgroup of $\Pic(\E)$ generated
by the components of the singular fibers.
\end{coro} 

\begin{lemm}\label{lemm:rank}
Let $\E\ra \P^1$ be a Jacobian elliptic fibration with 
singular fibers of simple multiplicative type (irreducible nodal curves). 
Let $M_1,M_2,M_3\subset S$ be three different irreducible 
horizontal divisors which are linearly independent 
in $\Pic(\E)$. Then at most one of the divisors
$M_{ij},$ ($i,j = 1,2,3,$ and $i\neq j$) is torsion.
\end{lemm}

{\em Proof.}
Indeed if two of the above divisors, for example $M_{12}$ and $M_{13}$,
are torsion then there is a positive integer
$N$ which annihilates both 
$\eta^*([M_{12}])$ and $ \eta^*([M_{13}])$.
Hence the kernel of $\eta^*$ 
has rank $\ge 2$. By assumption, the singular fibers
of $\E$ generate a subgroup of rank 1. Contradiction.   
\hfill $\square$

\

Let $\E\ra \P^1$ be a nonisotrivial elliptic fibration and 
$r$ the order of $\E$ in the Tate-Shafarevich
group of the corresponding Jacobian elliptic fibration 
${\cal J}(\E)$. 
Recall that for each integer $r'$ 
we have a principal homogeneous fibration  
${\cal J}^{(r')}(\E)$ of (relative) 
zero cycles of degree $r'$, together
with natural maps 
$$
{\cal J}^{(r')}(\E)\times {\cal J}^{(r'')}(\E)\ra {\cal J}^{(r'+r'')}(\E)
$$
and an identification of ${\cal J}^{(0)}(\E)$ and 
${\cal J}^{(r)}(\E)$, depending on the choice of a global section
of the Jacobian fibration ${\cal J}(\E)={\cal J}^{(0)}(\E)$.
After fixing the identification, we get
for every integer $m\equiv 1 \mod r$ 
a rational map 
$$
\phi_{m}\,:\, \E\ra \E
$$
of degree $m^2$, well defined up to the action of $H^0(\P^1,\E)$.  
This map is regular \'etale on the open (grouplike) part
$\E^0$ (complement to the singular points of the
singular fibers), but highly nonregular on $\E$.

\

\begin{lemm}\label{lemm:irreducible}
Let $\E\ra \P^1$ be a nonisotrivial elliptic fibration 
and $M$ a horizontal irreducible divisor. 
Then for all but finitely many primes (congruent to 
1 mod $r$) the preimage 
$\phi_p^{-1}(M)$ is also irreducible.
\end{lemm}

{\em Proof.} 
After base change, we can assume that $M$ is a nonzero section.
The global monodromy group of $\E_M$ over $M$ is a 
subgroup of a finite index in $\SL(2,\Z)$. 
Let ${\mathbb F}$ be the fundamental 
group of the complement $M \setminus {\rm Sing}$ where ${\rm Sing}$
is a subset of $M$ corresponding to singular fibers of $\E\ra M$.
Every section $M'\in \E_M$
defines a cocycle $s_{M'}\in H^1({\mathbb F}, (\Q/\Z)^2)$
(corresponding to the principal $(\Q/\Z)^2$-fibration, 
whose fiber over $b\in M\setminus {\rm Sing}$ 
is the set of points differing from $M'\cap \E_b$ by
torsion). 
This cocycle gives an affine action of ${\mathbb F}$ on $(\Q/\Z)^2$ 
and the orbits of the action of ${\mathbb F}$ 
correspond to the irreducible components of the preimages of $M'$ under
the maps $\phi_{m,M} : \E_M\to \E_M$ for all $m$.

If the section 
$M$ is not divisible in the group of sections $H^0(M,\E_M)$
by a prime $p$ and the monodromy map
$ {\mathbb F} \to \SL(2,\Z/p)$ is surjective then 
the preimage of $M$ under
the map $\phi_{p,M} : \E_M \ra \E_M$ is irreducible.
Indeed in this case either the action is affine and its orbit is 
$(\Z/p)^2$ or the action is linear.
In the first case the preimage of $M$ is irreducible. 
In the second case there is a section $M'$ with 
$ pM' = M$ and  hence $M$ is divisible as a section 
(see also \cite{bogomolov-tschinkel}). \hfill $\square$

\

\begin{lemm}\label{lemm:n1n2}
Let $\E\ra C$ be a nonisotrivial 
Jacobian elliptic fibration over an affine
connected and smooth curve $C$. 
Let $M_1, M_2$ be two irreducible multisections on $\E$ with the
property that there exists a smooth fiber $\E_b$  and a pair of
points $p_j\in \E_b\cap M_j$ such that $p_1-p_2$ is nontorsion in 
$\E_b$. Then there exists a positive integer
$N$ (depending on $\E,M_1,M_2$) such that
for all positive integers $n_1,n_2$ with $n_1+n_2>N$ the preimage 
$\phi_{n_1}^{-1}(M_1)$ intersects $\phi_{n_2}^{-1}(M_2)$.
\end{lemm}

{\em Proof.} 
The proof runs in the analytic category. 
It uses the following universal construction.

\

\noindent
{\bf Construction.} Consider the universal Jacobian elliptic curve
$u\,:\, \E_{\mathcal U}\ra {\mathcal H}$ which is obtained as a quotient 
${\mathcal H}\times \C/(1,\lambda)$, where $\lambda $ is a 
coordinate function in the upper halfplane ${\mathcal H}$.
This fibration is topologically trivial and has a natural trivialization
map   
$$
\kappa\,:\, {\mathcal H}\times \C/(1,\lambda)\ra {\rm T}=\C/(1,i),
$$ 
(where $i=\sqrt{-1}$). The map $\kappa$ is not complex analytic; however, 
the preimages of points $t\in {\rm T}$ are analytic sections of
the fibration 
$u\,:\, \E_{\mathcal U}\ra {\mathcal H}$ (this gives a nonanalytic family
of analytic sections). Indeed, the preimage of the point $t=a+bi$
is a section of $u$ which is given in a parametric 
form $(\lambda, b\lambda +a)\in {\mathcal H}\times \C$ (which descends to
$\E_{\mathcal U}$). 
The torsion sections of $u$ map into points in ${\rm T}$
(since $\kappa$ is a continuous homomorphism of algebraic groups).
The construction is equivariant with respect to $\SL(2,\Z)$.
Indeed, the action of $\SL(2,\Z)$ on $\E_{\mathcal U}$ transforms
torsion sections into torsion sections, (which are dense
in the family $\kappa^{-1}({\rm T})$). 
For any subgroup $\Gamma\subset \SL(2,\Z)$ 
of finite index which doesn't contain the center $\Z/2$
we have the induced map 
$$
u_{\Gamma}\,:\, \E_{\mathcal U}/\Gamma \ra {\mathcal H}/\Gamma
$$
and a factorization map $\pi_{\Gamma}$. 
We have the diagram:

\

\centerline{
\xymatrix{
\E_{\mathcal U}/\Gamma\ar[d]_{u_{\Gamma}}  
& \E_{\mathcal U} \ar[l]_{\pi_{\Gamma}} \ar[d]^{u}\ar[r]^{\kappa} & {\rm T} \\
{\mathcal H}/\Gamma  &  {\mathcal H} \ar[l]_{\pi_{\Gamma}}      &  
}
} 

\

\noindent
For every $t\in {\rm T}$, which is 
not a torsion point in ${\rm T}$ 
the orbit $\cup_{\gamma\in \Gamma} \gamma(t)$ is dense (in the
usual topology) in ${\rm T}$. Thus if $t$ is not torsion in ${\rm T}$
the intersection of the 
set $\pi_{\Gamma}(\kappa^{-1}(\cup_{\gamma\in \Gamma} \gamma(t)))$
with every fiber $\E_b\in \E_{\mathcal U}/\Gamma$ is dense. 

\

Now we return to the proof of Lemma~\ref{lemm:n1n2}. 
Consider the nonisotrivial Jacobian elliptic fibration
$\E\ra C$. 
We have a diagram

\

\centerline{
\xymatrix{
\E\ar[d]  & \E' \ar[l] \ar[d] & \E_U'\ar[l]\ar[d]\ar[r] & 
\E_{\mathcal U}\ar[d]\ar[r]& {\rm T} \\
C         &  C' \ar[l]  & C'\times_{{\mathcal H}/\Gamma}{\mathcal H}\ar[r]
\ar[l]& {\mathcal H} & 
}
} 

\

Changing the base, we reduce to 
the case when $M_1,M_2$ are sections of
a (topologically trivial) elliptic fibration 
$\E'\ra C'$ over some
analytic curve $C'$. 
Moreover, we can assume that the fibration $\E'\ra C'$
is induced from $\E_{\mathcal U}/\Gamma \ra {\mathcal H}/\Gamma$
under a dominant map $C'\ra {\mathcal H}/\Gamma$. 
We can identify $M_1'$ (the pullback of $M_1$) with
the zero section of $\E'$ 
(changing the zero section
amounts to changing the argument of $\kappa$ by a fiberwise translation). 
Then the image 
$M_{2,U}$ of $M_2$
in $\E_{\mathcal U}/\Gamma$ is algebraic, and therefore
not an orbit of $\Gamma$ on 
(the preimage of) some nontorsion point $t\in {\rm T}$. 
This means that
$\kappa(\pi_{\Gamma}^{-1}(M_{2,U}))$ 
covers some open (in the usual topology)
subset $V$ of ${\rm T}$. If $n_1+n_2$ is sufficiently big, then 
the translations of the set  $V$ by $n_2$ torsion points in ${\rm T}$ 
contain $n_1$ torsion points. \hfill $\square$

\subsection{Kummer surfaces}
\label{sect:kummer}

Let  $A$ be an abelian surface, 
$\tau$ the standard involution,
$A^*$ the blowup of $A$ in the 16 fixed points of $\tau$,  
$S$ the associated Kummer K3 surface and 
$\Ex$ the union of the 16 exceptional $(-2)$-curves
on $S$.  
We will say that  $A$ is {\em generic} if 
$\Pic(A)/\Pic^0(A) = \Z$ and if the endomorphism 
ring ${\rm End}(A) = \Z$.
In particular,  $A$ is not isogenous to
a product of elliptic curves.

\begin{lemm}\label{lemm:cond}
We can choose $A$ such that
\begin{enumerate}
\item $A$ is generic;
\item the elliptic fibration $S\ra \P^1$ is Jacobian;
\item all singular fibers of the elliptic fibration $S\ra \P^1$
are irreducible (consequently, 
the exceptional divisor $\Ex$ is horizontal in $S\ra\P^1$).  
\end{enumerate}
\end{lemm}

{\em Proof.}
The Picard lattice of a polarized Kummer surface  $S$ is
given by a sublattice (of rank at least 17) in 
$3H\oplus (-E_8)\oplus(-E_8)$, where $H$ is the standard
hyperbolic lattice. 
Let $h_A$ be the generator (polarization) of the Neron-Severi group 
of our generic abelian surface $A$. This class is invariant
under the involution. It descends to the singular surface
$\tilde{S}$ and lifts to a class $h_S\in \Pic(S)$.  
For any even positive integer $2k$ there exists 
a generic $A$ such that $h_S^2 = 2k$.   
For generic $A$ the Picard group of $S$ is a direct sum
$\Z h_S\oplus \Pi$, where $\Pi$ sits in the exact
sequence
$$
0\ra (-2){\rm Id}_{16}\ra \Pi\ra \Lambda^2((\Z/2)^4)\ra 0.
$$
The last projection extends naturally to a projection 
$\Pic(S)\ra\Lambda^2((\Z/2)^4)$ ($h_S$ is mapped to zero). 
Given a lattice vector of square zero we can find a translate
of this vector (under involutions of $\Pic(S)$
with respect to $(-2)$-classes) representing an 
elliptic fibration (see \cite{nikulin}). 

Remark that though the lattice 
$\Pi$ is not unimodular, for any primitive
element $e\in \Pic(S)$ with a nontrivial
projection to $\Lambda^2((\Z/2)^4)$ we can find an element $x\in \Pi$
such that $(e,x)=1$. 
Let $e'\in \Pic(S)$ be a class obtained from $e$ by 
reflexion with respect to a $(-2)$-class. Then there still exists a 
class $x'\in\Pic(S)$ with $(e',x')=1$. 

Choose a class $e$ 
of square zero giving an elliptic fibration $f_e\,:\,S\ra \P^1$, 
such that the projection of $e$ to  $\Lambda^2((\Z/2)^4)=(\Z/2)^6$
is nonzero. 
The lattice $N_e:=\{n\in \Pi\,|\, (n,e)=0\}$ 
is negative semi-definite.
We will choose $N_e$ such that it has no elements of
square $-2$.
This implies that the singular fibers of $f_e$ are
irreducible.
Simultaneously, we will choose
the (affine) lattice $L_e:=\{l\in\Pic(S)\,|\, (l,e)=1\}$ 
such that $L_e$ has $(-2)$-vectors.  
Since $N_e$ has no $(-2)$-classes  
every class (modulo translations by $e$) 
of square $(-2)$ in $L_e$  
corresponds to a section of $f_e\,:\, S\ra\P^1$.

The lattice $\Pi$ contains a finite number
of $(-2)$-vectors. Thus a generic vector in this lattice
is not orthogonal to any $(-2)$-vector. Take such a
primitive vector $x$   
and choose a polarization $h_S$ such that $h_S^2=-x^2$.
Now we can choose a (generic) abelian surface $A$
(with endomorphisms $\Z$)
such that the square of the generator $h_A$ of $\NS(A)$ equals $2h_S^2$.
On the corresponding Kummer K3 surface  $S$ we have 
$(h_S-x)^2=0$. It follows that 
$S$ has an elliptic fibration 
$f\,:\, S\ra \P^1$ and  that the lattice $N_f$
is isomorphic to $N_x\oplus \Z(h_S-x)$, and therefore
has no classes of square $(-2)$ (for any $x$ as above). 
Let $z$ be a vector $\Pi$ such  that $(x,z)=1$
(this is possible since $x$ is primitive). 
Every vector in $L_x$ is equal to $z+n + c(h_S-x)$
where $n\in N_x$ and $c\in \Z$. Its square is
$z^2+n^2 + 2(z,n) -2c$. Since the lattice $\Pi$ is even
and $c$ an arbitrary  integer it follows that $L_x$ always contains
classes of square $-2$. \hfill $\square$

\

\begin{rem}
The same proof shows that one can construct  Jacobian 
elliptic Kummer surfaces with all
singular fibers of simple multiplicative type without requiring 
that the associated abelian surface is generic. 
\end{rem}

We return to the general setup of Section~\ref{sect:construction}
and Section~\ref{sect:desingularization}.
Let $\E\ra \P^1$ be a Jacobian elliptic fibration, $M$ 
a horizontal divisor and $M_1,...,M_k$ its irreducible
components. Let $\Gamma(M)$ be the ``torsion'' graph
of $M$ -  each $M_l$ defines a vertex and 
two vertices $M_{l}$ and $M_{l'}'$
are connected by an edge 
if the difference divisor $M_{ll'}$ is not torsion.

\begin{lemm}\label{lemm:connected-ex}
Let $S$ be a Kummer surface as 
in the statement of Lemma~\ref{lemm:cond}.
Then the torsion graph $\Gamma(\Ex)$ is connected.
\end{lemm}

{\em Proof.} 
The divisor $\Ex$ consists of 16 irreducible  components,
their classes are linearly independent in $\Pic(S)$.
The rank of the subgroup of $\Pic(S)$ generated by 
the components of singular fibers is 1.
Applying Lemma~\ref{lemm:rank}
we see that every component of $\Ex$ is connected
to at least 14 other components of $\Ex$. 
In particular, every two vertices in $\Gamma(\Ex)$
are connected by a path of length at most 2. 
\hfill $\square$

\

\begin{lemm}\label{lemm:gamma-connected}
Consider $D_1^0,D_2^0\subset S^0$, obtained
as preimages of $\Ex$ under the maps $\phi_{q_1}, \phi_{q_2}$.
For almost all pairs of distinct prime numbers 
$q_1,q_2$ the intersection 
graph of the irreducible components of 
$D_1^0\cup D_2^0$ (in $S^0$) is connected. 
\end{lemm}

{\em Proof.}
By Lemma~\ref{lemm:irreducible},
the preimages of the irredicible components of $\Ex$ under
$\phi_{q_j}$ remain irreducible for almost all pairs of primes $q_1,q_2$.
Moreover, by Lemma~\ref{lemm:n1n2}, the {\em intersection} graphs
of the (irreducible) divisors $\phi^{-1}_{q_j}(\Ex_k)$ is 
connected for $q_j$ big enough. 
Now take an irreducible component $\Ex_1\subset\Ex$ and a component
$\Ex_k$, which is connected to $\Ex_1$ in the 
{\em torsion} graph $\Gamma(\Ex)$. 
Applying Lemma~\ref{lemm:n1n2} with $n_1=q_1$ and $n_2=q_2$ 
(and $q_1+q_2$ big enough) we see that
$\phi_{q_1}^{-1}(\Ex_1)$ intersects $\phi_{q_2}^{-1}(\Ex_k)$ (in $S^0$).
\hfill $\square$

\begin{rem}
The difficulty was to show that $D_1^0\cup D_2^0$ is
a connected divisor in the {\em open} surface $S^0$. 
The corresponding fact in the closed surface $S$ is trivial. 
\end{rem}

\subsection{Surfaces in abelian fourfolds}
\label{sect:morph}

We use the existence of nontrivial maps
between elliptic Kummer surfaces (see Sections~\ref{sect:fibr} 
and \ref{sect:kummer}) 
to construct interesting examples of surfaces
in abelian fourfolds. 

\

Fix a generic abelian surface $A$ such that
the associated Kummer surface $S$ admits a
Jacobian elliptic fibration 
(this is possible by Lemma~\ref{lemm:cond}).
Choose a pair $(q_1,q_2)$ of positive 
integers and consider the map
$$
\tilde{\phi}_{q_1,q_2}=(\tilde{\phi}_{q_1},\tilde{\phi}_{q_2}) 
\,:\, S\ra \tilde{S}\times \tilde{S}.
$$
Put $\tilde{g}_1=\tilde{\phi}_{q_1}, \tilde{g}_2=\tilde{\phi}_{q_2}$ and 
denote by $\tilde{X}=\tilde{X}_{q_1,q_2}$ and by $Y_1,Y_2,Y_{12}$ the surfaces 
obtained through the construction in 
Sections~\ref{sect:construction} and  \ref{sect:desingularization}. 
The surface $\tilde{X}$   
is lagrangian with respect to some 
nondegenerate 2-form on $A\times A$. 
We will show that for appropriate choices of 
integers $q_1,q_2$ the surface $X$ (and consequently $\tilde{X}$) 
is not fibered.

\subsection{Uniqueness}

\begin{prop}\label{prop:unique-form}
Let $A_1,A_2$ be simple abelian surfaces and 
$\tilde{X}$ a lagrangian surface in $A_1\times A_2$. 
Assume that $\tilde{X}$ is not isomorphic to an abelian surface,
that $\tilde{X}$ projects dominantly onto $A_1$ and $A_2$,  
that it is stable under 
the involutions $\tilde{\tau}_1,\tilde{\tau}_2$ and that
it is lagrangian with respect to at least two 
nonproportional nondegenerate 
forms, one of which is invariant with 
respect to both involutions $\tilde{\tau}_j$. 
Then $\tilde{X}$ is a finite unramified cover
of a product of two curves. 
\end{prop}

{\em Proof.}
By assumption, every $(1,0)$-form on $A_1\times A_2$ restricts
nontrivially to $\tilde{X}$. 
If there are two nonproportional nondegenerate  $(2,0)$-forms
which restrict trivially to $\tilde{X}$ then 
there exists a 
holomorphic 2-form $\omega$ of rank 2 on $A_1\times A_2$ 
which restricts trivially to $X$. 
(Indeed, the $(2,0)$-forms of rank 2 on the abelian variety
$A_1\times A_2$ correspond to points on a quadric in the projective
space $\P^5=\Proj(H^{2,0}(A_1\times A_2,\C))$. Any line in 
$\P^5$ intersects this quadric.)
Any such $\omega $ is equal to 
$w'\wedge w''$, where $w',w''$ are nonproportional
$(1,0)$-forms on $A_1\times A_2$.
Thus we have a pencil of curves ${\mathcal P}_{\omega}$ on $\tilde{X}$ 
such that both $w'$ and $w''$
are equal to zero on the fibers of this pencil. 
Therefore, 
we have a family of abelian surfaces $A_{t}\subset A_1\times A_2 $
(where $t$ is a point in the surface $B_{\omega}=(A_1\times A_2)/A_0$)
such that the intersection of $\tilde{X}$ with $A_t$ is either empty or
a fiber of the pencil ${\mathcal P}_{\omega}$.    
Hence, the image of $\tilde{X}$ under the projection to the base of 
${\mathcal P}_{\omega}$ coincides with the image of $\tilde{X}$ in
the abelian surface $B$. The image of $\tilde{X}$ in $B$ is a 
curve  $C_{\omega}$ of genus $\ge 2$. 
The form $\omega $ is induced from the holomorphic 
volume form on the abelian surface $B$. 
Consider the action of $\tilde{\tau}_1,\tilde{\tau}_2$ on $\omega$. 
Since this action transforms
$(2,0)$-forms 
of rank 2 on $A_1\times A_2$ onto themselves, we have the following
possibilities: either $\omega$ is invariant under $\Z/2+\Z/2$ 
(modulo multiplication by a constant)  or there 
is another 2-form  $\tilde{\omega}$ of rank 2, 
which is trivial on $\tilde{X}$. 
The first case is excluded by the assumptions
(that both projections on $A_1,A_2$ are dominant). 
In the second case 
we have another projection of $\tilde{X}$ onto an abelian surface, 
with fibers transversal to the fibers of the first projection.
Thus we have an isogeny 
$$
A_1\times A_2\ra  B_{\omega}\times B_{\tilde{\omega}},
$$
which exhibits $\tilde{X}$ as a finite
abelian covering of a product 
$C_{\omega}\times C_{\tilde{\omega}}$ of two curves of genus $\ge 2$. 
\hfill $\square$

\begin{coro}\label{coro:unique}
Let $\tilde{X}$ be a  surface obtained as
a $\Z/2+\Z/2$-cover of a (singular) K3-surface as 
in Section~\ref{sect:construction} (in particular, we do not assume 
that $S$ is elliptic). 
Assume that the conditions of Proposition~\ref{prop:unique-form}
hold for $\tilde{X}$. Then $\tilde{X}$ is isomorphic to $C\times C$, where
$C$ is a curve of genus 2 and the quotient map 
$$
\tilde{X}\ra S=(C\times C)/(\Z/2 + \Z/2)
$$
is described in Example~\ref{exam:unique}.
\end{coro}

{\em Proof.} 
Let $\tilde{X}$ be an unramified abelian cover of the product $C_1\times C_2$
of curves of genus $\ge 2$ as in the proof of 
Proposition~\ref{prop:unique-form}. 
Our assumption implies that there 
is a unique $(\Z/2 + \Z/2)$-equivariant $(2,0)$-form on $\tilde{X}$.
The involutions $\tau_1$ and $\tau_2$ interchange 
the projections of $X$ to $C_1\times C_2$ (since 
we have two projections of $X$ onto abelian surfaces $B_1,B_2$ and 
$A_1,A_2$ both map surjectively onto $B_1,B_2$ and 
the involution on $A_1$ interchanges the two projections). 
In particular, $C_1=C_2$ and $C_1\times C_1=C\times C$. 
The involution $\tau_{12}=\tau_1\tau_2$
induces an involution on the abelian cover of $C\times C$. 
Therefore, $\tau_{12}=(\sigma,\sigma)$, where $\sigma $ is some involution
on $C$. Thus we have a map 
$$
S=\tilde{X}/ (\Z/2 + \Z/2)\ra \Sym^2(C)/\sigma.
$$
The condition $h^{1,0}(S)=0$ implies that $\sigma $ is a hyperelliptic
involution and the condition $h^{2,0}(S)=1$ implies that $g(C)=2$.
Moreover, the map $\tilde{X}\ra C\times C$ is
in fact an isomorphism (unramified covers increase the Euler characteristic
and the number of invariant 2-forms has to increase as well). 
\hfill $\square$

\begin{coro}\label{coro:prod}
We keep the notations of Section~\ref{sect:construction}. 
Consider 
$$
\tilde{\varphi}=(\tilde{g}_1,\tilde{g}_2)\,:\, 
S\ra \tilde{S}_1\times \tilde{S}_2
$$
(and the associated maps $g_j\,:\, S\ra S_j$).  
Assume that at least one of the maps ${g}_1,{g}_2$ is not
an isomorphism or (if it is) the map $g_1g_2^{-1}$ doesn't 
lift to an automorphism of $C\times C$.  
Then $X$ is lagrangian 
with respect to exactly one (up to multiplication
by constants) nondegenerate form on $A_1\times A_2$.
\end{coro}

{\em Proof.} Indeed, the proof of Corollary~\ref{coro:unique}
shows that in the opposite case, $X$ is isomorphic to $C\times C$, 
with $g(C)=2$ and the action of $\Z/2 + \Z/2$ is the product
of the hyperelliptic involution and interchanging of coordinates, modulo
automorphisms of $C\times C$. \hfill $\square$
 
\

\subsection{1-forms}
\label{sect:1-forms}

\begin{lemm}\label{lemm:22}
Let $\hat{S}$ be birational to a K3 surface and 
$Y\ra \hat{S}$ its double cover. Assume that
$Y$ has at most $A1$-singularities (double points).
Then either $h^0(Y, \Omega^1)\le 2$ or
$\alb(Y)\subset\Alb(Y)$ is a hyperelliptic curve $C$ 
and the (covering) involution on $Y$ transforms into 
a hyperelliptic involution on $C$. 
\end{lemm}

{\em Proof.} Indeed, 
since there are no 1-forms on $S$, the involution acts  
by multiplication with $-1$ on $H^0(Y,\Omega^1)$. 
Thus the product of any two forms is invariant 
under the involution. 
Therefore, the map  
$$
\alpha\,:\, \wedge^2 H^0(Y,\Omega^1)\ra H^0(Y, \Omega^2)
$$
has image of dimension 1 or 0. If it is 0, then 
$\alb(Y)$ has dimension 1. (Indeed, it means that
any $(2,0)$-form on $\Alb(Y)$ pulls back to a $0$-form on $Y$
and hence $\alb(Y)$ has dimension 1.)
In this case the involution acts on $C=\alb(Y)\subset \Alb(Y)$
as well as on $H^0(C,\Omega^1)$ as $(-1)$. It follows that
$C$ is hyperelliptic and that the map $Y\ra C$ transforms
the involution on $Y$ into a hyperelliptic involution on $C$.

If the dimension of the image of $\alpha$ 
is 1, then there is a pair of nonproportional forms 
$\omega, \omega'$ such their product is nonzero. 
Assume that there is another form $\omega''$, linearly
independent of $\omega,\omega'$. Then $\omega'\wedge\omega''=
f\omega\wedge\omega'$, where $f$ is a nonconstant function, 
contradiction. \hfill $\square$

\begin{coro}
Let $Y_j$ be the double covers of $\hat{S}$ as above. Then
$h^0(Y_j,\Omega^1)=2$ (for $j=1,2$). Moreover, 
$Y_j$ admits no dominant rational maps onto
curves of genus $> 0$.
\end{coro}

{\em Proof.}
Indeed, both surfaces $Y_j$ 
admit a dominant map onto an abelian surface $A$, which has
two linearly independent
1-forms. Now we apply Lemma~\ref{lemm:22}. This proves the 
first statement. 
By assumption, the abelian surface $A$ contains
no elliptic curves. By the previous lemma, $\Alb(Y_j)$ is
isogenous to $A$. The second statement follows. \hfill $\square$. 

\

A priori, we don't know that $Y_{12}$ does not admit
dominant maps onto curves of genus $>0$. 
We have to consider the following (mutually exclusive) possibilities:

$h^0(Y_{12}, \Omega^1)=0$;

$h^0(Y_{12}, \Omega^1)=2$ and the wedge product of the two 
forms is nontrivial;

There is a projection $Y_{12}\ra C$, where $C$ is
a curve of genus $\ge 1$. 

\

We will show that the surfaces $X_{q_1,q_2}$ constructed in 
Section~\ref{sect:morph} are not of the last two types
for almost all pairs of primes $q_1,q_2$. The surfaces of
the first type are nonfibered; and there are   
simple examples such surfaces:

\

\begin{lemm}\label{lemm:y12}
Assume that the map $g_1\,:\, S\ra S_1$
is an isomorphism. 
Then $h^0(Y_{12}, \Omega^1) = 0$. 
\end{lemm}

{\em Proof.}
The covering $Y_1$ is (birational to) the abelian
surface $A_1$. The variety $X$ is (birational to)  
a fiber product of $A_1$ and some surface 
$Y_2$ over $\hat{S}$.
Consider the product 
$$
H^0(Y_2,\Omega^1)\times H^0(Y_{12},\Omega^1)\ra H^0(X,\Omega^2).
$$
The image is at most 1-dimensional, since it consists of 
$\tau_1$-invariant forms and since there is exactly one such 
form which is induced from $Y_1=A_1$. 
The space $H^0(Y_2,\Omega^1)$ has dimension 2 and hence
for any nontrivial form $\omega_{12}$ in $H^0(Y_{12},\Omega^1)$
its product with some form $\omega_2\in H^0(Y_2,\Omega^1)$
is zero. This means that $\omega_2$ is trivial on 
the fibers of some pencil of curves on $X$. This pencil covers
(birationally) $A_2$. The image of a generic
fiber of the pencil must be an elliptic curve in $A_2$
(being the zero set of a 1-form on $A_2$). 
This leads to a contradiction, 
since $A_2$ doesn't contain elliptic curves by assumption.
Thus there are no holomorphic 1-forms on $Y_{12}$.
\hfill $\square$

\begin{coro}\label{coro:g1-iso}
Assume that $g_1$ is an isomorphism and
that either $g_2$ is not an isomorphism or (if it is)
$g_1g_2^{-1}\,:\, S\ra S$ is not an involution on $S$.
Then the associated surface $X$ is not fibered 
and $\tilde{X}$ is lagrangian with
respect to exactly one 2-form on $A\times A$. 
\end{coro}

{\em Proof.} By Corollary~\ref{coro:unique}, $\tilde{X}$ is lagrangian
with respect to exactly one nondegenerate 2-form, unless
$\tilde{X}=C\times C$, where $C$ is a curve of genus 2. 
Thus it suffices to prove that $X$ is not fibered. 
Notice that $h^0(X,\Omega^1)=4$.
By Lemma~\ref{lemm:y12}, the Albanese $\Alb(X)$ is isogenous
to $A_1\times A_2$. The image  $\alb(X)$ in the Albanese is also lagrangian
with respect to exactly one nondegenerate 2-form. If $X$ were fibered
over a curve of genus $>1$ there would be another (a degenerate) 2-form on 
$\Alb(X)$ which would be trivial on $\alb(X)$. 
This is a contradiction. Notice also that $X$ doesn't admit
any dominant maps onto elliptic curves since $\Alb(X)$ doesn't map onto
elliptic curves. \hfill $\square$

\

In the following  sections we
will give further (more complicated) 
examples of nonfibered lagrangian surfaces in
$A_1\times A_2$.

\

\subsection{Weakly lagrangian structures}

Let $\tilde{X}\subset A_1\times A_2$ 
be as in Section~\ref{sect:construction} and $X$ its
(partial) desingularization as in Section~\ref{sect:desingularization}.  
Consider the set $L^X$ of (possibly degenerate) holomorphic
$(2,0)$-forms  
on $\Alb(X)$ which restrict trivially to $\alb(X)$.
It contains the set 
of weakly lagrangian 
structures on $\alb(X)$.
Since the condition is linear, $L^X$ is a linear subspace 
of $H^{2,0}(\Alb(X),\C)$,
The $\Z/2 + \Z/2$-action on $X$ lifts to an 
action on $\Alb(X)$, leaving the space $L^X$ invariant.

From now on we assume that $\tilde{X}$ is not isomorphic to 
$C\times C$,
where $C$ is a curve of genus $2$ (since in this case
$\dim L^X=2$, see Proposition~\ref{prop:unique-form}).
Furthermore, we will assume that $h^0(Y_{12},\Omega^1) \neq 0$
(if $\tilde{X}\neq C\times C$
and $h^0(Y_{12},\Omega^1)= 0$ then $\dim L^X = 1$). 

We have to consider the following cases:
\begin{enumerate}
\item  $h^0(Y_{12},\Omega^1) \neq 0$ and the wedge product
is degenerate;
\item 
$h^0(Y_{12},\Omega^1) = 2$ 
and the wedge product  
$$
\wedge\,:\, H^0(Y_{12},\Omega^1)\times H^0(Y_{12},\Omega^1)\ra
H^0(Y_{12},\Omega^2) 
$$ 
is nondegenerate. 
\end{enumerate}

\begin{rem}
Notice however, that 
examples of surfaces $X$ such that 
$$
h^0(Y_{12},\Omega^1) \neq 0
$$ 
are somewhat pathological. In particular, almost all surfaces
$X_{q_1,q_2}$ are not
of this type. 
\end{rem}

\

The space $L^X$ is invariant under $\Z/2 + \Z/2$ and hence 
decomposes into the direct sum 
$$
L^X = L^X_0 \oplus L^X_1 \oplus L^X_2 \oplus L^X_{12},
$$ 
where $L^X_0$ stands
for a subspace of $\Z/2 + \Z/2$-invariant forms in $L^X$, the space
$L^X_j$ is the space of $\tau_j$-invariant forms and 
$L^X_{12}$ is the space of $\tau_{12}$ invariant forms.

This decomposition arises from the decomposition of
the space $H^{2,0}$ of $(2,0)$-forms on $\Alb(X)$ 
under the $\Z/2 +\Z/2$-action.
We will denote these spaces
using the
same indices. 
We know that  
$$
\begin{array}{ccc}
H^{1,0}(\Alb(X))_0    & = & 0,\\
H^{1,0}(\Alb(X))_1    & = & H^{1,0}(A_2), \\
H^{1,0}(\Alb(X))_2    & = & H^{1,0}(A_1), \\
H^{1,0}(\Alb(X))_{12} & = & H^{1,0}(Y_{12}).
\end{array}
$$
Thus by taking the exterior product we obtain
that 
$$
H^{2,0}_0 = 
\Lambda^2 H^{1,0}(\Alb(X))_1 +  \Lambda^2 H^{1,0}(\Alb(X))_2
+  \Lambda^2 H^{1,0}(\Alb(X))_{12}.
$$
We also have the decompositions
$$
\begin{array}{ccc}
H^{2,0}(\Alb(X))_1    & = & H^{1,0}(\Alb(X))_{12}\times H^{1,0}(A_2),\\
H^{2,0}(\Alb(X))_2     & = & H^{1,0}(\Alb(X))_{12}\times H^{1,0}(A_2), \\
H^{2,0}(\Alb(X))_{12}  & = &  H^{1,0}(A_1)\times H^{1,0}(A_2).
\end{array}
$$

\begin{prop} 
\label{prop:cases}
If  $h^0(Y_{12},\Omega^1) \neq 0$ and the wedge product
is degenerate then  
$$
\dim L^X = \dim \Lambda^2 H^0(Y_{12}, \Omega^1) + 1.
$$
\end{prop}

{\em Proof.} We subdivide the proof into a sequence of lemmas. 

\begin{lemm} 
Under the conditions in Proposition~\ref{prop:cases} we have
$$
L^X_1 = L^X_2=L^X_{12} = 0.
$$
\end{lemm}

{\em Proof.} 
Since the space $L^X_{12}$ is induced from $A_1\times A_2$ 
it has dimension 0 (by Proposition~\ref{prop:unique-form} there is 
a unique lagrangian structure on $\tilde{X}\subset A_1\times A_2$).
We have
$$
L^X_1\subset H^{2,0}(\Alb(X))_1 = H^{1,0}(\Alb(X))_{12}\times 
H^{1,0}(A_2).
$$ 
Since  $\dim H^{1,0}(A_2) = 2$ for any element
$l\in L^X_1 $  there are two forms $w_1, w_2 \in H^{1,0}(A_2)$
and two forms $u_1,u_2 \in H^{1,0}(\Alb(X))_{12}$ with 
the property 
$$
w_1\wedge u_1 + w_2\wedge u_2 = l.
$$
We know that the forms $u_j$ are induced from a 
map $f : X\to C$, where $C$ is some hyperelliptic curve.
 
If $w_1, w_2$ are linearly dependent then 
$w_1$ is trivial on the  images of the fibers of the projection 
$f : X\ra C$. Since $X$ surjects upon $A_1$ and 
$A_2$ this means that $w_1$ is trivial on a
family of curves in $A_2$. This can happen only 
if $A_2$ is isogenous to a product of elliptic curves.
This contradicts the assumption that $A_2$ is simple.
 
Consider the case when $w_1,w_2$ are linearly
independent.
Since the wedge product on $H^{0,1}(\Alb(X))_{12}$ is trivial
we can write $u_2 = fu_1$ for some rational function $f$ on $X$ and
hence $ w_1 + fw_2\wedge u_1 = 0$. The function $f$
is constant on the fibers of the family $ X \ra C$ and the fibers
map into $A_1$. Since the generic fiber is not 
an elliptic curve the forms
$w_1,w_2\in H^{1,0}(A_2)$ are not linearly dependent 
on the fibers which
contradicts the assumption.

Similar argument yields $L^X_2 = 0$.
\hfill $\square$

\begin{coro}
Under the conditions of Proposition~\ref{prop:cases} we have 
$$
L^X = L^X_0= \Lambda^2 H^{1,0}(\Alb(X))_{12} + L_A,
$$
where
$\dim L_A = 1$. The subspace $L_A = 1$ is induced from $A_1\times A_2$. 
A generic element
in $L^X$ is nondegenerate if $\dim \Alb(X)$ is even and
is degenerate of corank $1$ if $\dim \Alb(X)$ is odd.
\end{coro}

{\em Proof.}
If 
$$
v = v_1 + v_2,\,\,  v_1\in L_A, \,\,v_2 \in 
\Lambda^2 H^{1,0}(\Alb(X))_{12}
$$ 
then
$\rk\, v = \rk\, v_1 + \rk\, v_2$ (since the spaces 
$ L_A, L^X_0$ belong to the second exterior powers
of complementary linear subspaces in $H^{1,0}(\Alb(X))$).

The nonzero element of $L_A$ is induced from a nondegenerate
form on $A_1\times A_2$ and has rank $4$. 
We are in the case when $Y_{12}$ (and consequently $X$)
admits a dominant map onto a curve $C$ of genus $\ge 1$.
If $g(C)$ is even the generic
element of $L^X$ is a nondegenerate form on $\Alb(X)$.
The dimension of $\Alb(X)$ is even iff $g(C)$ is even.
Similarly when $g(C)$ is odd.
\hfill $\square$

This finishes the proof of Proposition~\ref{prop:cases}.

\begin{lemm}\label{lemm:factor}
Under the conditions of Proposition~\ref{prop:cases} 
any dominant map $ h\,:\, X\ra C' $ of $X$ onto
a curve with $g(C')>0$ 
is a composition  
$$
h = sf,\, f\, :\, X\ra C,\,\,  s\,: \, C\ra C'.
$$
\end{lemm}

{\em Proof.}
Repeating the previous argument 
we see that 
the rank of a nonzero form in $L_A$ is $2$ and hence $\rk\, v \geq 2$
for 
$v = v_1 + v_2,$
$v_1\in L_A, $ $ v_2 \in \Lambda^2 H^{1,0}(\Alb(X))_{12} $
with $v_1\neq 0$. Thus 
any form of rank $1$ in $L^X$ belongs to 
$\Lambda^2 H^{1,0}(\Alb(X))_{12}$.
Since any map $h \,:\, X\ra C',\,\, g(C') \geq 2$ 
gives a nontrivial
form of rank $1$ we obtain the result in this case.

If $g(C') = 1$ then there is a map $\Alb(X) \to C'$.
The variety $\Alb(X)$ is isogenous to the product
$ \Jac(C)\times A_1\times A_2$ and any map from
$A_1\times A_2$ to a curve of genus $1$ is 
trivial. Hence the map $ \Jac(C)\times A_1\times A_2 \to C'$
is induced from the projection onto $\Jac(C)$.
It means that the holomorphic $(1,0)$-form lifted from $C'$
on $X$ is induced from $f : X \to C$.
This yields the result for $g(C') = 1$.
\hfill $\square$

\begin{coro} 
\label{coro:factoringc}
Under the conditions of Proposition~\ref{prop:cases}
every map $X\ra C'$ factors through the map 
$X\ra Y_{12}\ra C\ra C'$. 
\end{coro} 

\

\noindent
We consider the case when $h^0(Y_{12},\Omega^1) = 2$ 
and the wedge product on 1-forms is nondegenerate on $Y_{12}$.
It follows that the Albanese variety $\Alb(X)$ is isogenous
to a product of three abelian surfaces
to $A_1\times A_2\times A_{12}$. Indeed, for all 3 quotients $Y_1,Y_2,Y_{12}$
of $X$ by the nontrivial involutions the corresponding Albanese varieties
have dimension 2. These are very strong conditions; though we are not
sure that such examples exist, we 
would like to analyze this potential possibility:

\

In this case $\bar{Y}_{12}$ admits a 
$\Z/2$-equivariant map $\bar{Y}_{12}\ra A_{12}$, 
which descends to a map $\bar{S}\ra \tilde{S}_{12}$ (here 
$\tilde{S}_{12}$ is the K3 surface obtained as the quotient $A_{12}/\Z/2$).
The ramification over $S^0$ is the union of $D_1^0\cup D_2^0$. 
It follows that the preimage $D_{12}^0$ of $\Ex_{12}$ in $S^0$ 
(under the map $S^0\ra S_{12}$) is equal to $D_1^0\cup D_2^0$. 
It follows that $D_{12}^0$ is a smooth divisor and hence
$D_1^0$ and $D_2^0$ have no intersection in $S^0$. 

\begin{coro}
\label{coro:disconnect}
The intersection graph of the irreducible components of $D_{12}^0$
{\em in $S^0$} is totally disconnected. 
\end{coro}

\subsection{Nonfibered lagrangian surfaces}

In this section we show that for infinitely many pairs
of integers  $q_1,q_2$
the surface $X=X_{q_1,q_2}$ is not fibered 
(in particular, the lagrangian surface 
$\tilde{X}_{q_1,q_2}\subset A_1\times A_2$ is also not fibered).

\begin{lemm}\label{lemm:12}
Assume that $Y_{12}$ admits a dominant 
$\Z/2$-equivariant (with respect to the covering $Y_{12}\ra \hat{S}$) 
map $\rho$ onto a hyperelliptic curve $C$ of genus $>0$. 
Then $\rho$ descends to a dominant 
map $\bar{\rho}\,:\, \bar{Y}_{12}\ra C$
which is also equivariant with respect to the covering involution
$\bar{Y}_{12}\ra \bar{S}$. 
\end{lemm}

{\em Proof.} Indeed, the map $\rho$ is a composition of a map 
$Y_{12}\ra \Alb(Y_{12})$ and $\Alb(Y_{12})\ra \Jac(C)$. 
Now we apply Corollary~\ref{coro:factors}.\hfill $\square$

\

After factorization of $\bar{Y}_{12}$ by the involution $\bar{\tau}$
we obtain a map 
$\bar{\rho}_{\sigma}\,:\, \bar{S}\ra \P^1=C/\sigma$, 
(where $\sigma$ is the
hyperelliptic involution on $C$). Notice that the fibers
of $\bar{\rho}_{\sigma}$ are connected (since the 
same property holds for $\rho$). Denote by $R=\{p_1,...,p_{2g+2}\}\subset C$
the ramification divisor of $C\ra \P^1$.

\ 

\centerline{
\xymatrix{
 \bar{Y}_{12}\ar[d]_{\bar{\tau}}\ar[r]^{\bar{\rho}} & C\ar[d]^{\sigma}   \\
 \bar{S}\ar[r]_{\bar{\rho}_{\sigma}}   & \P^1   
} 
}

\

\begin{lemm}
The irreducible components 
of $D_1^0\cup D_2^0$ are contained in the fibers of $\bar{\rho}_{\sigma}$.
Moreover, they lie over $R\subset \P^1$.  
\end{lemm}

{\em Proof.} 
Since the map $\bar{\rho}$ commutes with the covering involution on
$\bar{Y}_{12}$ 
the image of the ramification divisor in $S^0\subset \bar{S}$ 
under the map $\bar{\rho}_{\sigma}$ 
is contained in $R$. 
The divisors $D_j^0$ are exactly the ramification divisors of
the double cover on $S^0$. 
\hfill $\square$

\begin{coro}\label{coro:connected}
Let $I_1^0\subset D_1^0$ and $I_2^0\subset D_2^0$ be 
two irreducible components of the divisors $D_1^0, D_2^0$. 
Assume that 
the intersection
$I_1^0\cap I_2^0\cap S^0\neq \emptyset$. 
Then $I_1^0$ and $I_2^0$ are mapped to the same point
in $R\subset \P^1$. 
\end{coro}

Now we describe the structure of those fibers of $\bar{\rho}_{\sigma}$ 
which lie over a point $p\in R$. 

\begin{lemm}\label{lemm:fibers}
All the components of a 
fiber in $S^0$ over a point $p\in R$, apart from components in $D_j^0$
have even multiplicities. 
\end{lemm}

{\em Proof.} Indeed, the double cover $\bar{Y}_{12}\ra \bar{S}$ is
induced from the double cover $C\ra \P^1$. Therefore, the ramification 
divisor of the former is the preimage of $R$. 
Since we know that the ramification divisor in $S^0$ is $D_1^0\cup D_2^0$
all the other components of the preimage of $R$ have even multiplicities.
\hfill $\square$

\begin{lemm}\label{lemm:y12-dominates}
Assume that the graph
of components of $D_1^0\cup D_2^0$ is connected
in $S^0$. Then there are no dominant maps 
of $Y_{12}$ onto curves of genus $>0$. 
Moreover, $\Alb(Y_{12})=0$.  
\end{lemm}

{\em Proof.}
Assume first that $\dim\Alb(Y_{12})=2$ and 
that the wedge product of 1-forms on $Y_{12}$ 
is nondegenerate. 
Then, by the Corollary~\ref{coro:disconnect} - this
contradicts the assumptions. Thus we have to consider
the remaining case when 
$Y_{12}$ does admit a dominant 
$\Z/2$-equivariant (with respect to the covering) 
map $\rho$ onto a hyperelliptic curve $C$ of genus $>0$.
Then $S^0$ has a regular map onto $\P^1$ such that
components of $D_1$ and $D_2$ are mapped onto 
ramification points $p_1,...,p_{2g+2}\subset R$. 
By Corollary~\ref{coro:connected},
all connected (in $S^0$) components of the 
union of the divisors $D_1^0$
and $D_2^0$ are mapped to one point, which we can assume to be $p_1$. 
Consider the singular symmetric tensor 
on $C$, which is equal to 0 at $\P^1$ and has the 
singularity $dz^{2n}/z^n$
at $p_2,...,p_{2g+2}$. 
Since the preimages of points $p_2,...,p_{2g+2}$ consist of components 
of multiplicity two this tensor lifts
into a  nonsingular symmetric tensor on $S^0$. 
Since the complement to $S^0$ in the K3 surface $S$ 
consists of finitely many points we can extend this tensor to
a holomorphic symmetric tensor on  $S$. However, there
are no such tensors on a K3 surface, contradiction. 
\hfill $\square$

\subsection{Maps to curves} 
\label{sect:maps} 

\begin{thm}\label{thm:main}
For almost all pairs of distinct prime numbers
$(q_1,q_2)$ the surface $X=X_{q_1,q_2}$ 
is not fibered.
\end{thm}

{\em Proof.}
By Lemma~\ref{lemm:gamma-connected}, 
we can insure (for almost all pairs
of distinct primes $q_1,q_2$) that the intersection
graph $D_1^0\cup D_2^0$ is connected in $S^0$. 
By Lemma~\ref{lemm:y12-dominates}, we have $\Alb(Y_{12})=0$. 
This implies that $X$ has a dominant map onto a curve only 
if $X$ is isomorphic to $C\times C$, where $C$ is a curve of genus 2. 
Contradiction. 
\hfill $\square$

\

\subsection{Fundamental groups}

\begin{prop}
Let $\tilde{X}\subset A_1\times A_2$ be 
a surface constructed in Section~\ref{sect:construction}.
Assume that $\tilde{X}$ is not isomorphic to $C\times C$, 
where $C$ is a curve of genus $\ge 2$. 
Then $\pi_1(\tilde{X})$ admits a surjective 
homomorphism onto a group $G$, where $G$ is a nontrivial central
extension of $\Z^8$ by $\Z^5$.
If $\Alb(\tilde{X})$ is isogenous to $A_1\times A_2$ then 
$$
G= \pi_1(\tilde{X})/
[[\pi_1(\tilde{X}),\pi_1(\tilde{X})] \pi_1(\tilde{X})]
$$ 
(modulo torsion). 
\end{prop}

{\em Proof.} For generic $A_j$ the $\Q$-linear envelope of
$\omega_j\in H^{2,0}(A_j)$ has dimension 5. Same holds
for $\lambda_2\omega_1 -\lambda_1\omega_2$ 
on $A_1\times A_2$ (for $\lambda_j\in\Q$). Thus the kernel is $\Z^5$. 
The image is $\Z^8=H^1(\tilde{X})/{\rm torsion}$.
It remains to recall a general fact from group cohomology:
Consider the exact sequences

\centerline{ 
\xymatrix{ 1\ar[r] &  K_1\ar[r] & H^2(G/[G,G],\Z)\ar[r] & H^2(G,\Z)\\
           1\ar[r] &  K_2\ar[r] & H^2(G/[G,G],\Z)\ar[r] & H^2(G/[[G,G],G],\Z)
}
}

Then $K_1$ and $K_2$ are isomorphic and, by duality,  
$K_2$ is isomorphic to $[G,G]/[[G,G],G]$ (tensor $\Q$).

\begin{rem}
There are two other potential cases:
when $Y_{12}$ admits a surjective map onto a curve $C$ of genus $g(C)=g$
and when $\Alb(\tilde{X})$ is isogenous to a product of 3 abelian surfaces.  
In the first case we have an exact sequence
$$
0\ra \Z^{1+i +g(2g-1)}\ra G\ra \Z^{8+2g}\ra 0.
$$
In the second case
$$
0\ra \Z^{4+2i}\ra G\ra \Z^{12}\ra 0.
$$
Here in both cases $0\le i\le 3$.
We believe that these cases do occur if $A_1,A_2$ have
Picard groups of higher rank or admit nontrivial endomorphisms.
Then the rational envelope of $(2,0)$-forms has smaller
dimension.  
\end{rem}

\end{document}